\documentclass[12pt]{amsart}
\usepackage{amsfonts}
\usepackage{amsmath,amssymb,amsthm}
\usepackage{color}

\newcommand{\R}{\mathbb R}

\newcommand{\C}{\mathbb C}

\newcommand{\ds}{\displaystyle}

\newtheorem{thm}{Theorem}[section]

\newtheorem{prop}[thm]{Proposition}
\theoremstyle{definition}

\theoremstyle{remark}

\begin{document}

\title[GENERAL ROTATIONAL SURFACES  IN $\R^4$]
{GENERAL  ROTATIONAL SURFACES  IN $\R^4$ WITH  MERIDIANS \\LYING
IN  TWO-DIMENSIONAL PLANES}

\author{Velichka Milousheva}

\address{Bulgarian Academy of Sciences, Institute of Mathematics and Informatics,
Acad. G. Bonchev Str. bl. 8, 1113, Sofia, Bulgaria; "L. Karavelov"
Civil Engineering Higher School, 175 Suhodolska Str.,  Sofia,
Bulgaria} \email{vmil@math.bas.bg}



\maketitle

\begin{abstract}
We apply the invariant theory of surfaces in the four-dimensional
Euclidean space to the class of general rotational surfaces with
meridians lying in two-dimensional planes. We find all minimal
super-conformal surfaces of this class.
\end{abstract}

\section{Introduction}
An invariant theory of surfaces in the four-dimensional Euclidean
space $\R^4$ was developed in \cite{GM1} and \cite{GM3} on the
base of the Weingarten map similarly to the classical theory of
surfaces in $\R^3$.

Let $M^2$ be a surface in $\R^4$ with tangent space $T_pM^2$ at
any point $p \in M^2$. In \cite{GM1} we introduced an invariant
linear map $\gamma$ of Weingarten-type at any $T_pM^2$, which
plays a similar role in the theory of surfaces in $\R^4$ as the
Weingarten map in the theory of surfaces in $\R^3$. This map
$\gamma$ generates two invariant functions $k$ and $\varkappa$,
analogous to the Gauss curvature and the mean curvature in $\R^3$.
The sign of  $k$ is a geometric invariant and the sign of
$\varkappa$ is invariant under the motions in $\R^4$. However, the
sign of $\varkappa$ changes under symmetries with respect to a
hyperplane in $\R^4$. The invariant $\varkappa$ is the curvature
of the normal connection of the surface $M^2$ in $\R^4$.

As in the classical case, the invariants $k$ and $\varkappa$
divide the points of $M^2$ into four types: flat, elliptic,
parabolic and hyperbolic. The surfaces consisting of flat points
satisfy the conditions $k=0, \,\, \varkappa=0$. These surfaces are
either planar surfaces ($M^2$  lies in a hyperplane $\R^3$ in
$\R^4$) or developable ruled surfaces in $\R^4$ \cite{GM1}.

The map $\gamma$ generates the second fundamental form at any
point $p \in M^2$. The first fundamental form and the second
fundamental form generate principal tangents and principal lines,
as in $\R^3$. Using the principal tangents, we obtained a
geometrically determined moving frame field at each point $p \in
M^2$. Writing derivative formulas of Frenet-type for this frame
field, we found eight invariant functions. In \cite{GM3} we proved
a fundamental theorem of Bonnet-type for surfaces without minimal
points, stating that these eight invariants determine the surface
up to a motion.

The minimal surfaces in $\R^4$ are characterized by the equality
$\varkappa^2 - k = 0.$ We proved in \cite{GM2} that on any minimal
surface $M^2$ the Gauss curvature $K$ and  the normal curvature
$\varkappa$ satisfy the following inequality $K^2-\varkappa^2\geq
0.$ This inequality generates two geometric classes of minimal
surfaces:
\begin{itemize}
\item the class of minimal super-conformal surfaces characterized by
$K^2 - \varkappa^2 =0$; \item the class of minimal surfaces of
general type characterized by $K^2-\varkappa^2>0$.
\end{itemize}

The class of minimal super-conformal surfaces in $\R^4$ is locally
equivalent to the class of holomorphic curves in $\C^2 \equiv
\R^4$ \cite{Ein}.

A fundamental theorem of Bonnet-type for strongly regular minimal
surfaces of general type was proved in \cite{GM2} in terms of four
invariant functions.

In the present paper we apply the invariant theory of surfaces on
a special class of surfaces in $\R^4$, which are general
rotational surfaces in the sense of C. Moore.

\section{Preliminaries}

Let $M^2: z = z(u,v), \, \, (u,v) \in {\mathcal D}$ (${\mathcal D}
\subset \R^2$) be a surface in $\R^4$ with tangent space $T_pM^2 =
{\rm span} \{z_u, z_v\}$ at an arbitrary point $p \in M^2$. We
choose an orthonormal normal frame field $\{e_1, e_2\}$ of $M^2$
so that the quadruple $\{z_u, z_v, e_1, e_2\}$ is positive
oriented in $\R^4$. Then the following derivative formulas hold:
$$\begin{array}{l}
\vspace{2mm} \nabla'_{z_u}z_u=z_{uu} = \Gamma_{11}^1 \, z_u +
\Gamma_{11}^2 \, z_v
+ c_{11}^1\, e_1 + c_{11}^2\, e_2,\\
\vspace{2mm} \nabla'_{z_u}z_v=z_{uv} = \Gamma_{12}^1 \, z_u +
\Gamma_{12}^2 \, z_v
+ c_{12}^1\, e_1 + c_{12}^2\, e_2,\\
\vspace{2mm} \nabla'_{z_v}z_v=z_{vv} = \Gamma_{22}^1 \, z_u +
\Gamma_{22}^2 \, z_v
+ c_{22}^1\, e_1 + c_{22}^2\, e_2,\\
\end{array}$$
where $\Gamma_{ij}^k$ are the Christoffel's symbols and $c_{ij}^k,
\,\, i,j,k = 1,2$ are functions on $M^2$.

Let $g$ be the standard metric in $\R^4$ and $\nabla'$ its flat
Levi-Civita connection.  We use the standard denotations
\;$E=g(z_u,z_u), \; F=g(z_u,z_v), \; G=g(z_v,z_v)$ for the
coefficients of the first fundamental form and set
$W=\sqrt{EG-F^2}$. Denoting by $\sigma$ the second fundamental
tensor of $M^2$, we have
$$\begin{array}{l}
\sigma(z_u,z_u)=c_{11}^1\, e_1 + c_{11}^2\, e_2,\\
[2mm]
\sigma(z_u,z_v)=c_{12}^1\, e_1 + c_{12}^2\, e_2,\\
[2mm] \sigma(z_v,z_v)=c_{22}^1\, e_1 + c_{22}^2\,
e_2.\end{array}$$

The three pairs of normal vectors $\{\sigma(z_u,z_u),
\sigma(z_u,z_v)\}$, $\{\sigma(z_u,z_u), \sigma(z_v,z_v)\}$,
$\{\sigma(z_u,z_v), \sigma(z_v,z_v)\}$ form three parallelograms
with oriented areas
$$\Delta_1 = \left|%
\begin{array}{cc}
\vspace{2mm}
  c_{11}^1 & c_{12}^1 \\
  c_{11}^2 & c_{12}^2 \\
\end{array}%
\right|, \quad
\Delta_2 = \left|%
\begin{array}{cc}
\vspace{2mm}
  c_{11}^1 & c_{22}^1 \\
  c_{11}^2 & c_{22}^2 \\
\end{array}%
\right|, \quad
\Delta_3 = \left|%
\begin{array}{cc}
\vspace{2mm}
  c_{12}^1 & c_{22}^1 \\
  c_{12}^2 & c_{22}^2 \\
\end{array}%
\right|,$$ respectively. These oriented areas determine three
functions
$$\ds{L = \frac{2 \Delta_1}{W}, \quad M = \frac{ \Delta_2}{W}, \quad N = \frac{2 \Delta_3}{W}},$$
which change in the same way as the coefficients $E, F, G$ under
any change of the parameters $(u,v)$.

Using the functions $L$, $M$, $N$ and $E$, $F$, $G$  we introduce
the linear map $\gamma$ in the tangent space at any point of $M^2$
$$\gamma: T_pM^2 \rightarrow T_pM^2$$
similarly to the theory of surfaces in $\R^3$.

The linear map $\gamma$ of Weingarten-type is invariant with
respect to changes of parameters on $M^2$ as well as to motions in
$\R^4$. Thus the functions
$$k = \frac{LN - M^2}{EG - F^2}, \qquad
\varkappa =\frac{EN+GL-2FM}{2(EG-F^2)}$$ are invariants of the
surface $M^2$.

The map $\gamma$ determines a second fundamental form of the
surface as follows. Let $X = \alpha z_u + \beta z_v, \,\,
(\alpha,\beta) \neq (0,0)$ be a tangent vector at a point $p \in
M^2$. The second fundamental form of $M^2$ at $p$ is defined by
$$II(\alpha,\beta) = - g(\gamma(X),X) = L\alpha^2 + 2M\alpha\beta + N\beta^2, \quad \alpha,\beta \in \R.$$
The notions of a normal curvature of a tangent, conjugate and
asymptotic tangents are introduced in the standard way by means of
$II$. The asymptotic tangents are characterized by zero normal
curvature.

A tangent $g: X = \alpha z_u + \beta z_v$ is said to be
\textit{principal} if it is perpendicular to its conjugate. A line
$c: u=u(q), \; v=v(q); \; q\in J \subset \R$ on $M^2$ is said to
be a \textit{principal line} if its tangent at any point is
principal. The surface $M^2$ is parameterized with respect to the
principal lines if and only if $F=0, \,\, M=0.$

Let $M^2$ be parameterized with respect to the principal lines and
denote the unit vector fields $\displaystyle{x=\frac{z_u}{\sqrt
E}, \; y=\frac{z_v}{\sqrt G}}$. The equality $M = 0$ implies that
the normal vector fields $\sigma(x,x)$ and $\sigma(y,y)$ are
collinear. We denote by $b$ a unit normal vector field collinear
with $\sigma(x,x)$ and $\sigma(y,y)$, and by $l$ the unit normal
vector field such that $\{x,y,b,l\}$ is a positive oriented
orthonormal frame field of $M^2$. Thus we obtain a geometrically
determined orthonormal frame field $\{x,y,b,l\}$ at each point $p
\in M^2$. With respect to the frame field $\{x,y,b,l\}$ we have
the following Frenet-type formulas:
$$\begin{array}{ll}
\vspace{2mm} \nabla'_xx=\quad \quad \quad \gamma_1\,y+\,\nu_1\,b;
& \qquad
\nabla'_xb=-\nu_1\,x-\lambda\,y\quad\quad \quad +\beta_1\,l;\\
\vspace{2mm} \nabla'_xy=-\gamma_1\,x\quad \quad \; + \; \lambda\,b
\; + \mu\,l;  & \qquad
\nabla'_yb=-\lambda\,x - \; \nu_2\,y\quad\quad \quad +\beta_2\,l;\\
\vspace{2mm} \nabla'_yx=\quad\quad \;-\gamma_2\,y \; + \lambda\,b
\; +\mu\,l;  & \qquad
\nabla'_xl= \quad \quad \quad \;-\mu\,y-\beta_1\,b;\\
\vspace{2mm} \nabla'_yy=\;\;\gamma_2\,x \quad\quad\quad+\nu_2\,b;
& \qquad \nabla'_yl=-\mu\,x \quad \quad \quad \;-\beta_2\,b,
\end{array}\leqno{(2.1)}$$
where  $\gamma_1, \gamma_2, \nu_1, \nu_2, \lambda, \mu, \beta_1,
\beta_2$  are geometric invariant functions.

The invariants $k$, $\varkappa$, and the Gauss curvature $K$ are
expressed by the functions $\nu_1$, $\nu_2$, $\lambda$, $\mu$ as
follows:
$$k = - 4 \nu_1 \nu_2 \mu^2, \qquad \varkappa = (\nu_1 - \nu_2) \mu, \qquad K = \nu_1 \nu_2 - (\lambda^2 + \mu^2). \leqno{(2.2)}$$

\vskip 5mm
\section{General rotational surfaces}

Considering general rotations in $\R^4$, C. Moore introduced
general rotational surfaces \cite{M} \, (see also \cite{MW1, MW2})
as follows. Let $c: x(u) = \left( x^1(u), x^2(u),  x^3(u),
x^4(u)\right); \,\, u \in J \subset \R$ be a smooth curve in
$\R^4$, and $\alpha$, $\beta$ are constants. A general rotation of
the meridian curve  $c$ in $\R^4$ is given by
$$X(u,v)= \left( X^1(u,v), X^2(u,v),  X^3(u,v), X^4(u,v)\right),$$
where
$$\begin{array}{ll}
\vspace{2mm} X^1(u,v) = x^1(u)\cos\alpha v - x^2(u)\sin\alpha v; &
\qquad
X^3(u,v) = x^3(u)\cos\beta v - x^4(u)\sin\beta v; \\
\vspace{2mm} X^2(u,v) = x^1(u)\sin\alpha v + x^2(u)\cos\alpha v;&
\qquad X^4(u,v) = x^3(u)\sin\beta v + x^4(u)\cos\beta v.
\end{array}$$
In the case $\beta = 0$ the $X^3X^4$-plane is fixed and one gets
the classical rotation about a fixed two-dimensional axis.

\vskip 2mm We consider a surface $\mathcal{M}^2$ in  $\R^4$,
defined by the vector-valued function
$$z(u,v) = \left( f(u) \cos\alpha v, f(u) \sin \alpha v, g(u) \cos \beta v, g(u) \sin \beta v \right);
\quad u \in J \subset \R, \,\,  v \in [0; 2\pi), \leqno{(3.1)}$$
where $f(u)$ and $g(u)$ are smooth functions, satisfying $\alpha^2
f^2(u)+ \beta^2 g^2(u) > 0$, $f'\,^2(u)+ g'\,^2(u) > 0$, $u \in
J,$ and $\alpha, \beta$ are positive constants. The surface
$\mathcal{M}^2$, given by (3.1) is a general rotational surface
whose meridians lie in  two-dimensional planes. In our case the
meridian is given by $m: x(u) = \left( f(u), 0,  g(u), 0\right);
\,\, u \in J \subset \R$.

Each parametric curve $u = u_0 = const$ of $\mathcal{M}^2$ is
given by
$$c_v: z(v) = \left( a \cos \alpha v, a \sin \alpha v, b \cos \beta v, b \sin \beta v \right);
\quad a = f(u_0), \,\, b = g(u_0)$$ and its Frenet curvatures are
$$\varkappa_{c_v} = \ds{\sqrt{\frac{a^2 \alpha^4 + b^2 \beta^4}{a^2 \alpha^2 + b^2 \beta^2}}}; \quad
\tau_{c_v} = \ds{\frac{ab \alpha \beta (\alpha^2 - \beta^2)}
{\sqrt{a^2 \alpha^4 + b^2 \beta^4}\sqrt{a^2 \alpha^2 + b^2
\beta^2}}}; \quad \sigma_{c_v} = \ds{\frac{\alpha \beta \sqrt{a^2
\alpha^2 + b^2 \beta^2}} {\sqrt{a^2 \alpha^4 + b^2 \beta^4}}}.$$
Hence, in case of $\alpha \neq \beta$ each parametric $v$-line is
a curve in $\R^4$ with constant curvatures (helix in $\R^4$
\cite{CDV}), and in case of $\alpha = \beta$ each parametric
$v$-line is a circle. We shall consider the case $\alpha \neq
\beta$.

Each parametric curve $v = v_0 = const$ of $\mathcal{M}^2$ is
given by
$$c_u: z(u) = \left(\, A_1 f(u), A_2 f(u), B_1 g(u), B_2 g(u) \, \right),$$
where  $A_1 = \cos \alpha v_0, \, A_2 = \sin \alpha v_0, \,B_1 =
\cos \beta v_0, \,B_2 = \sin \beta v_0$. The Frenet curvatures of
$c_u$ are expressed as follows:
$$\varkappa_{c_u} = \ds{\frac{|g' f'' - f' g''|}{(\sqrt{f'\,^2 + g'\,^2})^3}}; \quad \tau_{c_u} = 0.$$
Obviously, $c_u$ is a plane curve with curvature $\varkappa_{c_u}
= \ds{\frac{|g' f'' - f' g''|} {(\sqrt{f'\,^2 + g'\,^2})^3}}$. So,
for each $v = const$ the parametric curves $c_u$ are congruent in
$\R^4$. These curves are the \textit{meridians} of
$\mathcal{M}^2$.
 We shall call the surface, defined by (3.1) in the case $\alpha \neq \beta$, a \emph{general rotational surface}.

\vskip 2mm Calculating the tangent vector fields $z_u$ and $z_v$
we find the coefficients of the first fundamental form: $E =
f'\,^2(u)+ g'\,^2(u)$;  $F = 0$;  $G = \alpha^2 f^2(u)+ \beta^2
g^2(u).$ We consider the following orthonormal tangent frame field
$$\begin{array}{l}
\vspace{2mm}
x = \ds{\frac{1}{\sqrt{f'\,^2 + g'\,^2}}\left(f' \cos \alpha v, f' \sin \alpha v, g' \cos \beta v, g' \sin \beta v \right)};\\
\vspace{2mm} y = \ds{\frac{1}{\sqrt{\alpha^2 f^2 + \beta^2
g^2}}\left( - \alpha f \sin \alpha v, \alpha f \cos \alpha v, -
\beta g \sin \beta v, \beta g \cos \beta v \right)};
\end{array}$$
and the following orthonormal normal frame field
$$\begin{array}{l}
\vspace{2mm}
n_1 = \ds{\frac{1}{\sqrt{f'\,^2 + g'\,^2}}\left(g' \cos \alpha v, g' \sin \alpha v, - f' \cos \beta v, - f' \sin \beta v \right)};\\
\vspace{2mm} n_2 = \ds{\frac{1}{\sqrt{\alpha^2 f^2 + \beta^2
g^2}}\left( - \beta g \sin \alpha v, \beta g \cos \alpha v, \alpha
f \sin \beta v, - \alpha f \cos \beta v \right)}.
\end{array}$$
It is easy to verify that $\{x, y, n_1, n_2\}$ is a positive
oriented orthonormal frame field in $\R^4$.

Finding the second partial derivatives $z_{uu}$, $z_{uv}$,
$z_{vv}$, we calculate the functions $c_{ij}^k, \,\, i,j,k = 1,2$:
$$\begin{array}{ll}
\vspace{2mm} c_{11}^1 = g(z_{uu}, n_1) = \ds{\frac{g' f'' - f'
g''}{\sqrt{f'\,^2 + g'\,^2}}};
\quad \quad & c_{11}^2 =g(z_{uu}, n_2) = 0;\\
\vspace{2mm} c_{12}^1 = g(z_{uv}, n_1) = 0;
\quad \quad & c_{12}^2 = g(z_{uv}, n_2) = \ds{\frac{\alpha \beta (g f' - f g')}{\sqrt{\alpha^2 f^2 + \beta^2 g^2}}};\\
\vspace{2mm} c_{22}^1 = g(z_{vv}, n_1) = \ds{\frac{\beta^2 g f' -
\alpha^2 f g'}{\sqrt{f'\,^2 + g'\,^2}}}; \quad \quad & c_{22}^2 =
g(z_{vv}, n_2) = 0.
\end{array} \leqno{(3.2)}$$
Therefore the coefficients $L$, $M$ and $N$ of the second
fundamental form of $\mathcal{M}^2$ are expressed as follows:
$$L = \ds{\frac{2 \alpha \beta (g f' - f g') (g' f'' - f' g'')}{(\alpha^2 f^2 + \beta^2 g^2) (f'\,^2 + g'\,^2)}}; \qquad M = 0; \qquad
N = \ds{\frac{- 2\alpha \beta (g f' - f g') (\beta^2 g f' -
\alpha^2 f g')}{(\alpha^2 f^2 + \beta^2 g^2) (f'\,^2 +
g'\,^2)}}.$$

Consequently, the invariants $k$, $\varkappa$ and $K$ of
$\mathcal{M}^2$ are:
$$k = \ds{\frac{- 4 \alpha^2 \beta^2 (g f' - f g')^2 (g' f'' - f' g'') (\beta^2 g f' - \alpha^2 f g')}{(\alpha^2 f^2 + \beta^2 g^2)^3 (f'\,^2 + g'\,^2)^3}};$$

$$\varkappa =  \ds{\frac{\alpha \beta (g f' - f g')}{(\alpha^2 f^2 + \beta^2 g^2)^2 (f'\,^2 + g'\,^2)^2} \,
[(\alpha^2 f^2 + \beta^2  g^2)(g' f'' - f' g'') -  (f'\,^2 +
g'\,^2) (\beta^2 g f' - \alpha^2 f g') ]};$$

$$K =  \ds{\frac{(\alpha^2 f^2 + \beta^2  g^2)(\beta^2 g f' - \alpha^2 f g')(g' f'' - f' g'') - \alpha^2 \beta^2 (f'\,^2 + g'\,^2) (g f' - f g')^2}{(\alpha^2 f^2 + \beta^2 g^2)^2 (f'\,^2 + g'\,^2)^2} \,}.$$

In \cite{GM3} we found all general rotational surfaces  consisting
of parabolic points, i.e.  $k = 0$.

\vskip 2mm Now we shall apply the invariant theory of surfaces in
$\R^4$ finding the geometric  invariant functions  $\gamma_1,
\gamma_2, \nu_1, \nu_2, \lambda, \mu, \beta_1, \beta_2$ in the
Frenet-type formulas of $\mathcal{M}^2$.

The positive oriented orthonormal frame field $\{x, y, n_1, n_2\}$
defined above is the geometric frame field of $\mathcal{M}^2$, the
$u$-lines and $v$-lines of $\mathcal{M}^2$ are principal lines.

Let $\sigma$ be the second fundamental tensor of $\mathcal{M}^2$.
Using (3.2) we obtain
$$\begin{array}{l}
\vspace{2mm}
\sigma(x, x) = \ds{\frac{g' f'' - f' g''}{\left(\sqrt{f'\,^2 + g'\,^2}\right)^3}\,\, n_1};\\
\vspace{2mm}
\sigma(x, y) = \ds{\frac{\alpha \beta (g f' - f g')}{\sqrt{f'\,^2 + g'\,^2}(\alpha^2 f^2 + \beta^2 g^2)}\,\, n_2};\\
\vspace{2mm} \sigma(y, y) = \ds{\frac{\beta^2 g f' - \alpha^2 f
g'}{\sqrt{f'\,^2 + g'\,^2}(\alpha^2 f^2 + \beta^2 g^2)}\,\, n_1}.
\end{array}$$
The partial derivatives of the normal vector field $n_1$ are:
$$\begin{array}{l}
\vspace{2mm}
(n_1)_u = \ds{\frac{f' g'' - g' f''}{\left(\sqrt{f'\,^2 + g'\,^2}\right)^3}\left(f' \cos \alpha v, f' \sin \alpha v, g' \cos \beta v, g' \sin \beta v \right)};\\
\vspace{2mm} (n_1)v = \ds{\frac{1}{\sqrt{f'\,^2 + g'\,^2}}\left( -
\alpha g' \sin \alpha v, \alpha g' \cos \alpha v, \beta f' \sin
\beta v, - \beta f' \cos \beta v \right)}.
\end{array}$$
The last equalities imply
$$\langle (n_1)_u, n_2 \rangle = 0; \qquad \langle (n_1)_v, n_2 \rangle = \ds{\frac{\alpha \beta
(f f' + g g')}{\sqrt{f'\,^2 + g'\,^2}\sqrt{\alpha^2 f^2 + \beta^2
g^2}}}.$$

So, the invariants in the Frenet-type derivative formulas (2.1)
for the surface $\mathcal{M}^2$ are given by
$$\begin{array}{ll}
\vspace{2mm}
\gamma_1 = 0; & \qquad \nu_1 = \ds{\frac{g' f'' - f' g''}{(f'\,^2 + g'\,^2)^{\frac{3}{2}}}};\\
\vspace{2mm} \gamma_2 = \ds{- \frac{\alpha^2 f f' + \beta^2 g
g'}{\sqrt{f'\,^2 + g'\,^2}(\alpha^2 f^2 + \beta^2 g^2)}}; & \qquad
\nu_2 = \ds{\frac{\beta^2 g f' - \alpha^2 f g'}{\sqrt{f'\,^2 + g'\,^2}(\alpha^2 f^2 + \beta^2 g^2)}};\\
\vspace{2mm}
\lambda = 0; & \qquad  \beta_1 = 0;\\
\vspace{2mm} \mu = \ds{\frac{\alpha \beta (g f' - f
g')}{\sqrt{f'\,^2 + g'\,^2}(\alpha^2 f^2 + \beta^2 g^2)}}; &
\qquad  \beta_2 = \ds{\frac{\alpha \beta (f f' + g
g')}{\sqrt{f'\,^2 + g'\,^2}\sqrt{\alpha^2 f^2 + \beta^2 g^2} }}.
\end{array} \leqno{(3.3)}$$

Thus we obtain that for the class of general rotational surfaces
only five invariants from the general theory are essential (they
may not be zero). These five invariant functions determine the
general rotational surfaces up to a motion in $\R^4$.

\vskip 5mm
\section{Minimal super-conformal general rotational surfaces}

In this section we shall find all minimal super-conformal general
rotational surfaces.

\vskip 2mm

Let us recall that the \textit{ellipse of normal curvature} at a
point $p$ of a surface $M^2$ in $\R^4$ is the ellipse in the
normal space at the point $p$ given by $\{\sigma(x,x): \, x \in
T_pM^2, \, g(x,x) = 1\}$ \cite{MW1, MW2}. Let $\{x,y\}$ be an
orthonormal base of the tangent space $T_pM^2$ at $p$. Then, for
any $v = \cos \psi \, x + \sin \psi \, y$, we have
$$\sigma(v, v) = H + \ds{\cos 2\psi  \, \frac{\sigma(x,x) - \sigma(y,y)}{2}
+ \sin 2 \psi  \, \sigma(x,y)},$$ where $H = \ds{\frac{\sigma(x,x)
+ \sigma(y,y)}{2}}$ \, is the mean curvature vector of $M^2$ at
$p$. So, when $v$ goes once around the unit tangent circle, the
vector $\sigma(v,v)$ goes twice around the ellipse centered at
$H$. Obviously, $M^2$ is minimal if and only if for each point $p
\in M^2$ the ellipse of curvature is centered at $p$.

A surface $M^2$ in $\R^4$ is called \textit{super-conformal}
\cite{BFLPP} if at any point of $M^2$ the ellipse of normal
curvature is a circle. An explicit construction of any simply
connected super-conformal surface in $\R^4$ that is free of
minimal and flat points was given in \cite{DT}.

 \vskip 2mm
In \cite{GM2} it was proved that the class of minimal
super-conformal surfaces is characterized by the equalities
$$ \varkappa^2 - k = 0; \qquad K^2 - \varkappa^2 = 0. \leqno{(4.1)}$$

Without loss of generality we assume that the meridian $m$ is
given by
$$\begin{array}{ll}
\vspace{1mm}
f = u;\\
\vspace{1mm} g = g(u).
\end{array}$$
Then, from (3.3) we obtain that the functions $\nu_1$, $\nu_2$,
$\mu$ are expressed by the function $g(u)$ and its derivatives as
follows:
$$\nu_1 = \ds{\frac{- g''}{(1 + g'\,^2)^{\frac{3}{2}}}}; \quad
\nu_2 = \ds{\frac{\beta^2 g  - \alpha^2 u g'}{\sqrt{1 +
g'\,^2}(\alpha^2 u^2 + \beta^2 g^2)}}; \quad
 \mu = \ds{\frac{\alpha \beta (g  - ug')}{\sqrt{1 + g'\,^2}(\alpha^2 u^2 + \beta^2 g^2)}}. \leqno{(4.2)}$$

Now, using (2.2), (4.1) and (4.2) we get that the general
rotational surface $\mathcal{M}^2$, given by (3.1), is minimal
super-conformal if and only if $g(u)$ satisfies the equation:
$$\alpha \beta (g - u g') = \varepsilon (\alpha^2 u g' - \beta^2 g), \qquad \varepsilon = \pm 1. \leqno{(4.3)}$$

All solutions of equation (4.3) are given by
$$g(u) = \ds{c\,u^{\varepsilon \frac{\beta}{\alpha}}},\qquad c = const.$$

Thus we obtained that the class of all minimal super-conformal
general rotational surfaces is described as follows:

\begin{prop}
Let $\mathcal{M}^2$ be a general rotational surface in $\R^4$.
Then $\mathcal{M}^2$ is  minimal super-conformal if and only if
the meridian is determined by
$$g(u) = c\,u^{\varepsilon \frac{\beta}{\alpha}},\qquad c = const.$$
\end{prop}

\vskip 3mm Denoting by $k = \varepsilon \ds{\frac{\beta}{\alpha}}$
\, ($k \neq \pm 1$), we obtain that all minimal super-conformal
general rotational surfaces are given by
$$\mathcal{M}^2: z(u,v) = \left( u \cos\alpha v,\, u \sin \alpha v, \,c\, u^k \cos \beta v,\, c\, u^k \sin \beta v \right).$$

The invariants $k$, $\varkappa$ and $K$ of $\mathcal{M}^2$ are
expressed as follows:
$$k = 4 \ds{\frac{c^4 k^4(1-k)^4 u^{4(k-2)}}{\left(1 + c^2 k^2\,u^{2(k-1)}\right)^6}}; \quad
\varkappa = 2\, \varepsilon  \ds{\frac{c^2 k^2(1-k)^2
u^{2(k-2)}}{\left(1 + c^2 k^2\,u^{2(k-1)}\right)^3}}; \quad K = -
2 \ds{\frac{c^2 k^2(1-k)^2 u^{2(k-2)}}{\left(1 + c^2
k^2\,u^{2(k-1)}\right)^3}}.$$

\end{document}